# Distributions that are both log-symmetric and R-symmetric

## M.C. Jones

*Department of Mathematics & Statistics, The Open University*
*Milton Keynes MK7 6AA, U.K.*
*e-mail:* m.c.jones@open.ac.uk

## Barry C. Arnold

*Department of Statistics, University of California – Riverside*
*Riverside CA 92521, U.S.A.*
*e-mail:* barry.arnold@ucr.edu

**Abstract:** Two concepts of symmetry for the distributions of positive random variables $Y$ are log-symmetry (symmetry of the distribution of $\log Y$) and R-symmetry [7]. In this paper, we characterise the distributions that have both properties, which we call doubly symmetric. It turns out that doubly symmetric distributions constitute a subset of those distributions that are moment-equivalent to the lognormal distribution. They include the lognormal, some members of the Berg/Askey class of distributions, and a number of others for which we give an explicit construction (based on work of A.J. Pakes) and note some properties; Stieltjes classes, however, are not doubly symmetric.

**AMS 2000 subject classifications:** Primary 62E10; secondary 60E05.

**Keywords and phrases:** double symmetry, lognormal distribution, moment equivalence, weighted distribution.

## 1. Introduction

Let $Y$ be a random variable following an absolutely continuous distribution $F$ with density $f$ on the positive half-line $\mathbb{R}^+$. One natural concept of symmetry of this distribution is that $Y/\delta$ and $\delta/Y$ have the same distribution for some $\delta > 0$. Since this is equivalent to ordinary symmetry, about $\log \delta$, of the distribution of $\log Y$, we will refer to this as log-symmetry ([6] Section 12.A.c, [11]). In density terms,

$$y^2 f(\delta y) = f\left(\frac{\delta}{y}\right), \quad \text{for some } \delta > 0 \text{ and all } y > 0. \qquad (1)$$

Note that $\delta$ is the *median* of the log-symmetric distribution $F$.

An alternative interesting concept of symmetry on $\mathbb{R}^+$ is the R-symmetry introduced and investigated by Mudholkar and Wang [7]. This is defined directly in terms of the density and another constant $\theta > 0$ by

$$f(\theta y) = f\left(\frac{\theta}{y}\right), \quad \text{for some } \theta > 0 \text{ and all } y > 0. \qquad (2)$$

In fact, $\theta$ is the *mode* of $F$ (if $f$ is unimodal). A comparison of these two concepts of reciprocal symmetry was given by Jones [5].

The focus of this note is on distributions on $\mathbb{R}^+$ which have both properties (1) and (2), log-symmetry and R-symmetry; we will call such distributions *doubly symmetric*. It is easy to see that no continuous distribution on $\mathbb{R}^+$ can be doubly symmetric with the same constant ($\delta = \theta$). However, as noted in [5], the lognormal distribution — with density

$$f_0(y) = \frac{1}{\sqrt{2\pi}\sigma y} \exp\left\{-\frac{1}{2}\frac{(\log y - \mu)^2}{\sigma^2}\right\} \tag{3}$$

in its usual parametrisation resulting from taking the exponential of a random variable following a normal distribution with mean $\mu$ and variance $\sigma^2$ — is both log-symmetric about $\delta_1 = e^\mu$ and R-symmetric about $\theta_1 = e^{\mu-\sigma^2}$. In that paper, the first author of this note was "tempt[ed] to conjecture that the lognormal is unique in this respect" — this proves to be wrong!

In this note, we characterise the class of absolutely continuous distributions on $\mathbb{R}^+$ which are doubly symmetric (call this class DS). It turns out to be a proper subset of the wide class of absolutely continuous distributions on $\mathbb{R}^+$ which are moment-equivalent to the lognormal distribution (e.g. [9], Section 3; use MEL to denote this class); moreover, this subset of MEL distributions contains more densities than just the lognormal. See Section 2 for our main result. Our work will rely heavily on deep and detailed work of Pakes [8, 9].

The first, and most famous, MEL subfamily is a version of the Stieltjes class with density

$$f_\epsilon(y) = f_0(y)\left[1 + \epsilon \sin\left\{2\pi(\log y - \mu)/\sigma^2\right\}\right], \tag{4}$$

$-1 \leq \epsilon \leq 1$. This is the class considered by Stieltjes [12] and Heyde [4]; see [13] for a more general formulation. It is intriguing to find that such Stieltjes densities are *not* DS. In fact, members of the Stieltjes class (4) are neither log-symmetric nor R-symmetric except when $\epsilon = 0$. Instead, $y^2 f_\epsilon(\delta_1 y) = f_{-\epsilon}(\delta_1/y)$ and $f_\epsilon(\theta_1 y) = f_{-\epsilon}(\theta_1/y)$, the former comprising Theorem 3.1 of [9]. And it is clear that no other values of $\delta$ or $\theta$ can be found to make $f_\epsilon$ log-symmetric or R-symmetric, respectively, for $\epsilon \neq 0$.

For examples of MEL densities that *are* DS, and some investigation of their properties, see Section 3.

The paper finishes in Section 4 with some further comments.

## 2. Main result

Write $k = \delta/\theta$. Note that $k > 1$ because the median is not less than the mode for either log-symmetric or R-symmetric densities [5]. Then, manipulating (1) and (2) shows that DS densities satisfy

$$\frac{1}{k^2}f\left(\frac{y}{k^2}\right) = \frac{y^2 f(y)}{\theta^2 k^4} \tag{5}$$

(together with either of (1) and (2)). This, therefore, says that DS densities are the same when you rescale by $k^2$ (left-hand side of (5)) and when you perform weighting with weight function $y^2$ (right-hand side of (5); note that $E(Y^2)$ is indeed equal to $\theta^2 k^4$).

Rescaling, let $g(x) = \theta k^2 f(\theta k^2 x)$ be the density of $X = Y/(\theta k^2)$. Then, (5) translates to

$$\frac{1}{k^2} g\left(\frac{y}{k^2}\right) = y^2 g(y) \tag{6}$$

which is the same as (5) except that $E(X^2) = 1$. Now, (6) is identical to (3.15) of [8] provided we set Pakes's $q = k^{-2}$. Next, write $h(z) = g(\sqrt{z})/(2\sqrt{z})$ to be the density of $Z = X^2$. Then, as Pakes notes, the problem of solving (6) for $g$ is equivalent to solving

$$\frac{1}{k^4} h\left(\frac{y}{k^4}\right) = yh(y) \tag{7}$$

with $E(Z) = 1$.

Pakes ([8], Theorem 3.1) gives the general form for any distribution function solving (7) with $E(Z) = 1$. These distributions, it turns out [10], are a subset of the distributions that are moment-equivalent to the log-normal law. We are interested only in absolutely continuous versions of these distributions, for which Pakes [8] writes down the form of the density, $h^*$ say, at his (3.10) (see also Theorem 3.3 of [9]). It has a piecewise form on intervals $(k^{-4(i+1)}, k^{-4i}]$, $i = 0, \pm 1, \pm 2, \ldots$. Explicitly,

$$h^*(y) \propto \sum_{i=-\infty}^{\infty} k^{2i(i+3)} y^i \omega(k^{4i} y) I(k^{-4(i+1)} < y \leq k^{-4i}) \tag{8}$$

where $\omega$ is a nonnegative function on $(k^{-4}, 1]$ such that $h^*$ is integrable. Since $f(y) = (2y/\theta^2 k^4) h\left(y^2/(\theta^2 k^4)\right)$, solutions of (5) take the form

$$f^*(y) \propto \sum_{i=-\infty}^{\infty} \theta^{-2i} k^{2i(i+1)} y^{2i+1} \omega(\theta^{-2} k^{4(i-1)} y^2) I(\theta k^{-2i} < y \leq \theta k^{2-2i}). \tag{9}$$

However, densities $f^*$ constitute a superset of DS densities since the latter also need to satisfy one of (1) or (2).

**Theorem 1.** *Absolutely continuous doubly symmetric densities have the form $f^*$ given by (9) with $\omega$ chosen to satisfy*

$$\psi(u) = \psi\left(\frac{1}{k^4 u}\right), \qquad k^{-4} < u \leq 1. \tag{10}$$

*Here,*

$$\psi(u) \equiv u\,\omega(u).$$

*Proof.* First, write (9) as

$$f^*(y) \propto \sum_{i=-\infty}^{\infty} \theta^{-2i} k^{2i(i+1)} y^{2i+1} \mathcal{F}(\theta^{-2} k^{4(i-1)} y^2)$$

where $\mathcal{F}(u) = \omega(u) I(k^{-4} < u \le 1)$. Now, consider requirement (5), recalling that $\delta = k\theta$. We have

$$y^2 f^*(k\theta y) \propto \sum_{i=-\infty}^{\infty} k^{2i(i+2)} y^{2i+3} \mathcal{F}(k^{4i-2} y^2). \tag{11}$$

Also,

$$\begin{aligned} f^*(k\theta/y) &\propto \sum_{j=-\infty}^{\infty} k^{2j(j+2)} y^{-2j-1} \mathcal{F}(k^{4j-2}/y^2) \\ &= \sum_{i=-\infty}^{\infty} k^{2i(i-2)} y^{2i-1} \mathcal{F}(k^{-4i-2}/y^2). \end{aligned} \tag{12}$$

Expression (12) follows from the previous line by setting $j = -i$, the choice determined to be the one that gives the $\mathcal{F}$ functions in (11) and (12) the same support. It follows that (11) = (12) if

$$k^{8i} y^4 \omega(k^{4i-2} y^2) = \omega(k^{-4i-2}/y^2)$$

which, by setting $u = k^{4i-2} y^2$, is equivalent to (10). The reader can verify that the same result is obtained by considering (6) in place of (5). □

*Remarks.* Note that the mode $\theta$ enters (9) simply as a scale parameter. Note also that, although $f^*$ is the density of the (scaled) square root of a random variable whose distribution is MEL, it remains the case that $f^*$ itself is MEL; this is because the (scaled) square root of a lognormal random variable is itself lognormal. It will also prove useful, in Section 3, to write

$$f^*(\theta y) \propto \sum_{i=-\infty}^{\infty} k^{2i(i-1)} y^{2i-1} \psi(k^{4(i-1)} y^2) I(k^{-2i} < y \le k^{2-2i}). \tag{13}$$

## 3. Examples of DS distributions

### 3.1. Lognormal distribution

Lognormal density $f_0$ given by (3) or any alternative scaling of it by factor $c$ say, has $E(Y^2) = c^2 e^{2\mu+2\sigma^2} = \theta^2 k^4$ since $\theta = ce^{\mu-\sigma^2}$ and $k = e^{\sigma^2}$. Note that the $i = 1$ term of (13) shows that, on support $(k^{-4}, 1]$, $\psi(u) \propto f^*(\theta\sqrt{u})/\sqrt{u}$. Then, the log-normal's psi-function is

$$\psi(u) \propto \frac{1}{\sqrt{u}} \exp\left(-\frac{1}{8} \frac{\log^2 u}{\log k}\right) \tag{14}$$

which satisfies (10).

### *3.2. Another known MEL distribution*

For a first example of a class of MEL distributions some of whose members are DS and are not the lognormal, we turn to the Askey/Berg class ([1, 2, 8] Section 3, especially pp. 834–835, [9] Section 3). One way of writing the density of this class, employing a particular scaling for convenience, is as

$$f_\gamma(y) \propto y^{\gamma-1}/L_k(y) \quad \text{where} \quad L_k(y) = \sum_{n=-\infty}^{\infty} y^n k^{-\frac{1}{2}n^2}, \qquad (15)$$

$\gamma \in \mathbb{R}$, $y > 0$, is Ramanujan's theta function (e.g. [3]). Now,

$$\begin{aligned}
\frac{L_k(k^c y)}{L_k(k^c/y)} &= \frac{\sum_{n=-\infty}^{\infty} y^n k^{-\frac{1}{2}n^2+cn}}{\sum_{n=-\infty}^{\infty} y^{-n} k^{-\frac{1}{2}n^2+cn}} = \frac{\sum_{n=-\infty}^{\infty} y^n k^{-\frac{1}{2}n^2+cn}}{\sum_{n=-\infty}^{\infty} y^n k^{-\frac{1}{2}n^2-cn}} \\
&= \frac{\sum_{n=-\infty}^{\infty} y^n k^{-\frac{1}{2}(n-c)^2}}{\sum_{n=-\infty}^{\infty} y^n k^{-\frac{1}{2}(n+c)^2}} = y^{2c} \frac{\sum^1 y^j k^{-\frac{1}{2}j^2}}{\sum^2 y^j k^{-\frac{1}{2}j^2}}.
\end{aligned}$$

Here, $\sum^1$ sums over $j \in \{\ldots, -2-c, -1-c, -c, 1-c, 2-c, \ldots\}$ and $\sum^2$ sums over $j \in \{\ldots, -2+c, -1+c, c, 1+c, 2+c, \ldots\}$. The two sums are equal if $c$ is integer or half-integer, in which case (1) holds with $\delta = k^\gamma$ and (2) with $\theta = k^{\gamma-1}$ if $\gamma$ is also integer or half-integer. This chimes with (1.10) of Theorem 1 of [3] whose integer $n$ is our $2c$. We have thus shown that Askey/Berg densities with integer or half-integer values of $\gamma$ are DS.

In fact, limited computational comparisons of the lognormal and Askey/Berg densities with the same median and mode (as well as moments) have failed to find any visible differences between them. (Often, the two densities are identical to the eighth decimal place, sometime to just four.) The version of the lognormal density that is compared with $f_\gamma$ at (15) has $\mu = \gamma \log k$, $\sigma^2 = \log k$ and hence density

$$f_0(y) = \frac{1}{\sqrt{2\pi \log k}} k^{-\frac{1}{2}\gamma^2} y^{\gamma-1} \exp\left(-\frac{1}{2}\frac{\log^2 y}{\log k}\right).$$

It follows that, to a very good degree of approximation, it must be that

$$L_k(y) \propto \exp\left(\frac{1}{2}\frac{\log^2 y}{\log k}\right).$$

Indeed, it can be proven that $f_\gamma$ and $f_0$ are exactly equal at all points of the form $k^p$, $p = 0, \pm 1, \pm 2, \ldots$. To see this, note that, exactly,

$$L_k(k^p) = \sum_{n=-\infty}^{\infty} k^{-\frac{1}{2}n^2+np} = e^{\frac{1}{2}p^2} \sum_{n=-\infty}^{\infty} k^{-\frac{1}{2}(n-p)^2} \propto e^{\frac{1}{2}p^2}.$$

### 3.3. Constructing further specific DS MEL distributions

Appropriate functions $\psi$ satisfying (10) can be constructed by specifying nonnegative $\psi$ arbitrarily on $[k^{-2}, 1]$ and extending it uniquely to the remainder of $[k^{-4}, 1]$ by defining

$$\psi(u) = \psi\left(\frac{1}{k^4 u}\right), \qquad k^{-4} \leq u \leq k^{-2}. \tag{16}$$

Note that $\psi$ is continuous at $k^{-2}$. Also, it can readily be shown from (9) that continuity of density $f^*$ at the "grid endpoints" $k^{2i}, i = 0, \pm 1 \pm 2, \ldots$ requires $\psi(1) = \psi(k^{-4})$ — equivalent to Lemma 3.1 of [9] — which is immediate from (10); it follows that if $\psi$ is chosen to be continuous then $f^*$ is continuous. We note Pakes's ([9], p. 1278) observation that continuous MEL densities "probably are the exception" and that this construction seems a promising way of extending this class (with double symmetry thrown in).

Now let $\psi$ be continuously differentiable. From (16), continuous differentiability of $f^*$ at "grid midpoints" $k^{2i-1}, i = 0, \pm 1, \pm 2, \ldots$ requires $\psi'(k^{-2}) = 0$. From (13), continuous differentiability of $f^*$ at grid endpoints $k^{2i}, i = 0, \pm 1, \pm 2, \ldots$ requires $2\psi'(1) = -\psi(1) < 0$.

Unimodality of $f^*$ with unique mode $\theta$ demands that:

(i) $u^{i-(1/2)}\psi(u)$ be increasing for all $i = 1, 2, \ldots$ which follows if $\sqrt{u}\psi(u)$ is strictly increasing;

(ii) $u^{i-(1/2)}\psi(u)$ be decreasing for all $i = 0, -1, -2, \ldots$ which follows if $\psi(u)/\sqrt{u}$ is strictly decreasing.

It can then be shown that, when $\psi$ is differentiable, these two requirements reduce to choosing $\psi$ on $(k^{-2}, 1)$ such that $|(\log \psi)'(u)| < (2u)^{-1}$.

One parametric family of unimodal DS (MEL) distributions arises if we take

$$\psi(u) = \psi_\alpha(u) \equiv u^{\alpha-(1/2)}, \qquad k^{-2} \leq u \leq 1, \qquad 0 < \alpha < 1. \tag{17}$$

We note at once that, while continuous, the resulting density is not continuously differentiable at both grid mid- and end-points for any value of $\alpha$. However, this at least makes for densities that more obviously differ from the lognormal. To this end, we briefly investigate the simplest case, that of $\alpha = 0$ in (17). The corresponding density is

$$f_U^*(y) = \frac{\sum_{i=-\infty}^{\infty} k^{2i(i-1)}(y/\theta)^{2i-1} I(k^{-2i} < y/\theta \leq k^{2-2i})}{\theta\{2\log k + \sum_{j=1}^{\infty} j^{-1} k^{-2j^2}(k^{2j} - k^{-2j})\}}. \tag{18}$$

It is pictured in Figures 1 and 2 in conjunction with the moment-equivalent lognormal density in each case. Density (18) acts as a piecewise (odd-power) polynomial approximation to the lognormal with discontinuities in derivative at grid endpoints. Discrepancies between the two densities are clearest at their common modes.

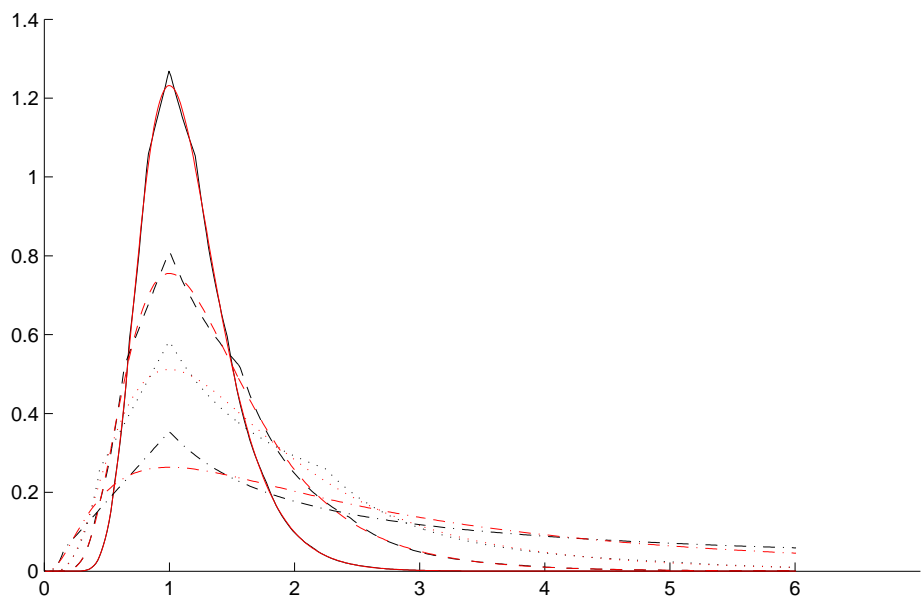

FIG 1. *Densities $f_U^*(y)$ (in black) and $f_0(y)$ (in red) for $\theta = 1$ and $k = 1.1$ (solid lines), $k = 1.25$ (dashed lines), $k = 1.5$ (dotted lines) and $k = 2.5$ (dot-dashed lines).*

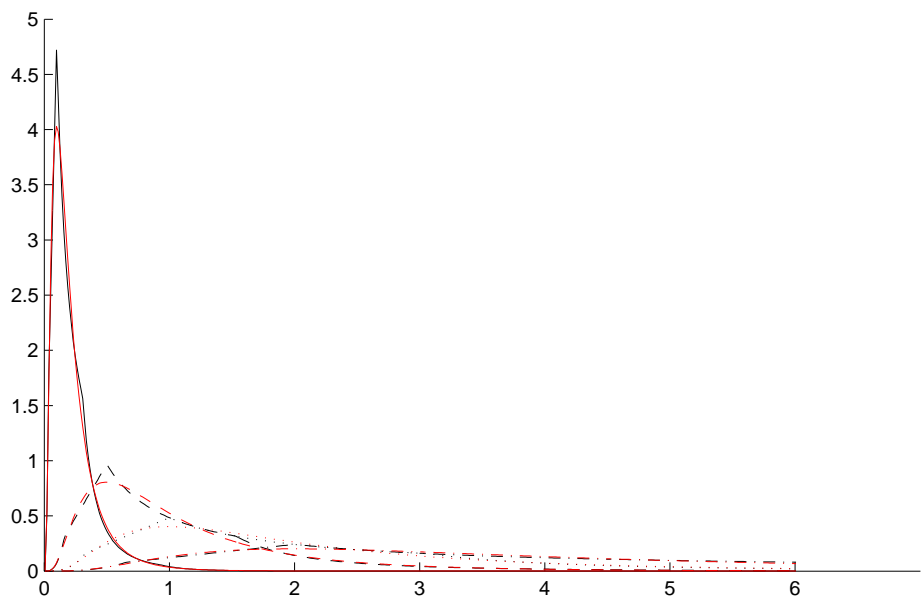

FIG 2. *Densities $f_U^*(y)$ (in black) and $f_0(y)$ (in red) for $k = 1.75$ and $\theta = 0.1$ (solid lines), $\theta = 0.5$ (dashed lines), $\theta = 1.0$ (dotted lines) and $\theta = 2.0$ (dot-dashed lines).*

## 4. Further remarks

That doubly symmetric distributions are moment-equivalent to the lognormal distribution is essentially immediate from the derived relationship (5) from which we get

$$E(Y^{s+2}) = \frac{\delta^{2(s+2)}}{\theta^{2(s+1)}} E(Y^s), \qquad s \in \mathbb{R}. \tag{19}$$

(Inter alia, the existence of all positive and negative moments of $f$ is thereby established.) This relationship is satisfied by the lognormal distribution with $\delta = \exp(\mu)$, $\theta = \exp(\mu - \sigma^2)$ and $E(Y^s) = \exp(s\mu + \frac{1}{2}s^2\sigma^2)$. In the closely related context of (7), Pakes [8, 9] notes that the ratio of the moment functions of $f$ and of the lognormal is periodic.

Here is another take on periodic differences between solutions. Define $h(w) = \log f(e^w), w \in \mathbb{R}$. Then, a little manipulation based on (5) yields

$$h(w - 2\log k) = 2w - 2\log\theta - 2\log k + h(w). \tag{20}$$

Thus, $h$ comprises a periodic function, with period $2\log k$, plus a linear one. Moreover, for DS densities we also have, from (2), that

$$h(w) = h(2\log\theta - w), \tag{21}$$

displaying symmetry of $h$ about $2\log\theta$. Combining (21) and (20) yields

$$\ell(w - 2\log k) = \ell(2\log\theta - w) \tag{22}$$

where $\ell(w) = h(w) - w$.

The periodicity evident in (20)–(22) can be utilized to formulate a characterization of the lognormal distribution within the class DS. If $X$ has a lognormal distribution then for any $\gamma > 0$, $X^\gamma$ is also a lognormal variable and consequently is also DS. If one assumes only that for any $\gamma > 0$, $X^\gamma \in DS$, it follows that the functions in (20)–(22) have multiple periodicities and, as a consequence, the function $h(w)$ must be quadratic. In this manner we may formulate:

**Theorem 2.** *If $X^\gamma \in DS \; \forall \gamma > 0$, then $X$ has a lognormal distribution.*

We conjecture, but are as yet unable to verify, that it will suffice to replace in the hypothesis the statement "$\forall \gamma > 0$" by the statement "for two judiciously chosen values of $\gamma$".

## References

[1] ASKEY, R. (1989). Orthogonal polynomials and theta functions. In *Theta Functions – Bowdoin 1987, Proceedings of Symposia in Pure Mathematics* **49** (eds: R.C. Gunning & L. Ehrenpreis), pp. 299–321; American Mathematical Society. MR1013179

[2] BERG, C. (1998). From discrete to absolutely continuous solutions of indeterminate moment problems. *Arab J. Math. Sci.* **4** 1–18. MR1667218

[3] BERNDT, B.C. (1993). Ramanujan's theory of theta-functions. In *Theta Functions; From the Classical to the Modern* (ed: M.R. Murty), pp. 1–64; American Mathematical Society. MR1224050

[4] HEYDE, C.C. (1963). On a property of the lognormal distribution. *J. Roy. Statist. Soc. Ser. B* **29** 392–393. MR0171336

[5] JONES, M.C. (2008). On reciprocal symmetry. *J. Statist. Planning Inference* **138** 3039–3043.

[6] MARSHALL, A.W. and OLKIN, I. (2007). *Life Distributions; Structure of Nonparametric, Semiparametric, and Parametric Families.* Wiley, Hoboken, NJ. MR2344835

[7] MUDHOLKAR, G.S. and WANG, H. (2007). IG-symmetry and R-symmetry: interrelations and applications to the inverse Gaussian theory. *J. Statist. Planning Inference* **137** 3655–3671. MR2363286

[8] PAKES, A.G. (1996). Length biasing and laws equivalent to the log-normal. *J. Math. Anal. Applic.* **197** 825–854. MR1373083

[9] PAKES, A.G. (2007) Structure of Stieltjes classes of moment-equivalent probability laws. *J. Math. Anal. Applic.* **326** 1268–1290. MR2280980

[10] PAKES, A.G. and KHATTREE, R. (1992). Length-biasing, characterization of laws, and the moment problem. *Austral. J. Statist.* **34** 307–322 MR1193781

[11] SESHADRI, V. (1965). On random variables which have the same distribution as their reciprocals. *Canad. Math. Bull.* **8** 819–924. MR0192530

[12] STIELTJES, T.J. (1894/5) Recherches sur les fractions continues. *Ann. Faculté Sci. Toulouse* **8** 1–122; **9** 5–47. MR1508159

[13] STOYANOV, J. (2004). Stieltjes classes for moment-indeterminate probability distributions. *J. Appl. Probab.* **41A** 281–294. MR2057580